\documentclass{amsart}
\usepackage{tabularx}
\usepackage{mathrsfs}
\usepackage{amsmath,amssymb,amsfonts,amsthm}
\usepackage{caption}
\usepackage{ctable}
\usepackage{geometry}
\newtheorem{theorem}{Theorem}

\title{Applying the symmetry groups to study the $n$ body problem}
\author{Zhihong Xia, Tingjie Zhou}
\address{Department of Mathematics, Northwestern University, Evanston,
  IL 60208 USA}
\address{Department of Mathematics, Southern University of Science and 
  Technology, Shenzhen, China}
\email{xia@math.northwestern.edu , 11930530@mail.sustech.edu.cn}
\date{November 4, 2021}

\begin{document}

\maketitle   

\begin{abstract}
  We introduce an algebraic method to study local stability in the
  Newtonian $n$-body problem when certain symmetries are present. We
  use representation theory of groups to simplify the calculations of
  certain eigenvalue problems. The method should be applicable in many
  cases, we give two main examples here: the square central
  configurations with four equal masses, and the equilateral
  triangular configurations with three equal masses plus an additional
  mass of arbitrary size at the center.  Using representation theory
  of finite groups, we explicitly found the
  eigenvalues of certain $8\times 8$ Hessians in these examples, with
  only some simple calculations of traces. We also studied the local
  stability properties of corresponding relative equilibria in the
  four-body problems.
   
\

   \textbf{Keywords:} the Newtonian $n$-body problem, symmetry, representation theory, eigenvalues

\end{abstract}


\section{Introduction}
In this article, we consider the Newtonian $n$-body problem with
certain symmetries. We study the local stability problems of certain
relative equilibria corresponding to symmetric central
configurations. The idea is actually very classical in physics
literature and we bring it to the $n$-body problem: the representation
theory of finite groups can be used effectively to simplify the
calculation of eigenvalue problems associated with central
configurations and relative equilibria. This technique provides us
with some new perspectives in celestial mechanics, and
helps us to obtain some results that is difficult to obtain
otherwise. There are many interesting scenarios where this technique
might apply, we start with two very specific examples.

A relative equilibrium in the $n$-body problem is where the $n$ point
masses maintain a fixed relative position during their motions. Only a
special set of configurations of the $n$ bodies can maintain such
relative equilibria, these are called {\em central
  configurations}\/. Central configurations are well-studied objects
(cf.\ \cite{2012Four, 2019On, 2017On,Palmore1976, Saari1980On,
  smale1998mathematical, xia2004convex}), dating back to Euler and
Lagrange, a simple question whether the number of central
configurations is finite, for any given set of $n$ positive masses,
remains a  major open problem in mathematics (cf.\ Smale
\cite{smale1998mathematical}, Albouy and Kaloshin\cite{2012Finiteness},
Hampton \& Moeckel \cite{HamptonMoeckel2006} ).

As an easy example, $n$ equal masses placed on vertices of a regular
$n$-gon form a central configuration. One can also place one more mass
of any size at the center of  the above regular $n$-gon, this will form
an central configuration for $n+1$ body problem.  In fact, J.\ C.\ Maxwell famously
studied this configuration as a model for Saturn's ring system in
1859. Moeckel
\cite{MR1350320} and Roberts \cite{MR1816912} provided detailed
analysis of corresponding relative equilibria and obtained
the linear stability results of these relative equilibria. Dihedral symmetry
played an important role in their analysis.

In this paper, we systematically use symmetry and group
representations to study some of these models. Specifically, we
consider two examples, with different symmetries, in the planar
four-body problem. The first example is a simple central configuration
where we place four equal masses on corners of a square. The second
example is where three masses are placed at vertices of an equilateral
triangle, and the fourth particle of arbitrary mass, placed at the
center. As is well-known, central configurations are critical points
of certain scaled potential function. To classify each of these
central configurations, one needs to compute the eigenvalues of
certain $8 \times 8$ matrices, typically not directly solvable
analytically. Moreover, the second example involves an arbitrary
parameter, the mass at the center, complicating any attempt at
analytical or numerical solutions. In both cases, we are able to give
a simple analytical solution involving only simple calculation of
traces, by using the representation theory of associate symmetry
groups, the dihedral groups of equilateral triangle and square.

Interestingly, Palmore\cite{Palmore1975Classifying, Palmore1976} showed
that a degeneracy arose in the $4$-body central configurations. The existence of
degeneracies in relative equilibria answers several questions raised
by Smale. In this work, by giving an analytical formula, we find exactly
where and how Palmore's degeneracy occurs.

We obtained the complete local stability analysis for the 
relative equilibrium corresponding to the square configuration. We
also discussed how symmetry and group representation can simplify
local stability analysis for the other case.

Leandro~\cite{leandro2017factorization} introduced the
representation theory to study factorization of characteristic
polynomial in various models including Maxwell's ring system. He applied
his method to a Rhombus family of relative equilibria
(Leandro~\cite{MR3951830}). What we present here is based on the method
introduced by Xia \cite{Xia2008}.

\

\noindent {\bf Acknowledgment:} The authors thank Yanxia Deng and Fang
Wang for their helpful 
  discussions. The authors are also grateful to the referee for many 
  useful suggestions.

\section{PRELIMINARIES}
In this section, we introduce some basic concepts of the central
configurations as well as the main ideas from the representation theory of
finite groups.

\subsection{Planar Central Configurations}
Let $q_i \in \mathbb{R}^2$ be the position of the particle $m_i$,
$i=1,2,\dots,n$. If there is a constant $\lambda$ such that
$$ \lambda m_iq_i = \sum_{1\le i<j \le n} \frac{m_im_j}{|q_i-q_j|^3}(q_j-q_i)$$
for all $1 \le i \le n$, then the $n$ particles are said to form a
(planar) {\em central configuration}\/. A planar central configuration
remains a central configuration after a rotation in $\mathbb{R}^2$ and
a scalar multiplication (cf.\ Xia \cite{xia2004convex}). Rotation and
scaling generate an {\em equivalent classes}\/ of central configurations.

It turns out that central configurations are critical points of the function $\sqrt{I}U$, where
$$I=\frac{1}{2}\sum_{i=1}^{n}m_i|q_i|^2,$$
$$U=\sum_{1 \le i<j \le n} \frac{m_im_j}{|q_i-q_j|}.$$
The function $I$ is the moment of inertia of the $n$-body system, while $U$
is the potential function. Equivalently, central configurations can be described as critical points of
the function U on the ellipsoid \cite{xia2004convex}
$$S=\{q=(q_1,\dots,q_n)\in R^{2n}|I=1\}.$$

\subsection{Hamiltonian of the $n$ Body Problem}
Let $p_i=m_i \dot{q_i}$, then the $n$-body problem has a Hamiltonian
formulation with the Hamiltonian~\cite{easton1993introduction}
$$ H=K-U, $$
where
$$K=\sum_{i=1}^{n} \frac{|p_i|^2}{2m_i}.$$
And the equations of motion are
$$\dot{q_i}=\frac{\partial H}{\partial p_i}=\frac{p_i}{m_i},$$
$$\dot{p_i}=-\frac{\partial H}{\partial q_i}=\frac{\partial
  U}{\partial q_i}.$$

There are some special solutions where the $n$ particles move on
concentric circles with uniform angular
velocity~\cite{1981Periodic}. Let the center of circle be at the
origin, then the solutions have the form 
$$q_i^*=exp(-\omega Jt)a_i,$$
$$p_i^*=-m_i\omega J exp(-\omega Jt)a_i,$$
where $a_i$ is constant vector, 
$$J=\begin{pmatrix}
0 & 1 \\
-1 & 0
\end{pmatrix}$$
and
$$exp(-\omega Jt)=\begin{pmatrix}
cos(\omega t) & -sin(\omega t) \\
sin(\omega t) & cos(\omega t)
\end{pmatrix},$$
where the angular velocity  $\omega$ is a positive constant that satisfies
\begin{equation}
\omega^2a_i+\sum_{\substack{j=1 \\ j\ne
    i}}^{n}\frac{m_j(a_j-a_i)}{|a_j-a_i|^3}=0. \label{omega}
\end{equation}
Solutions $a_1,\dots,a_n$ to the above equation, corresponding to positions for the
$n$ particles, are precisely the central
configurations.  Under a rotating coordinate system with constant
frequency $\omega$, these solutions will be fixed points, and therefore called {\em
  relative equilibria}\/.
Let $q_i=exp(-\omega Jt)x_i$ and $p_i=exp(-\omega Jt)y_i$. The
Hamiltonian, under uniform rotating coordinates, becomes
$$H=\sum_{i=1}^{n}\big( \frac{|y_i|^2}{2m_i}-\omega x_i^{T}Jy_i\big)
-\sum_{1 \le i<j \le n}\frac{m_im_j}{|x_i-x_j|}.$$
The equations of motion are
\begin{equation}
\begin{aligned}
\dot{x_i} &=\omega Jx_i +\frac{y_i}{m_i}, \\
\dot{y_i} &=\omega Jy_i+\sum_{\substack{j=1 \\ j \ne i}
}^{n}\frac{m_im_j(x_j-x_i)}{|x_j-x_i|^3}.
\end{aligned}
\end{equation}
The equations can be written as a second-order equation by eliminating $y_i$:
\begin{equation}
\ddot{x_i}=2\omega J\dot{x_i}+\omega^2
x_i+\frac{1}{m_i}\sum_{\substack{j=1 \\ j \ne i}
}^{n}\frac{m_im_j(x_j-x_i)}{|x_j-x_i|^3}.\label{con:inventoryflow} 
\end{equation}

\subsection{Groups and Symmetry}
Let $G$ be a group. The following are useful examples in our applications: 
\begin{itemize}
	\item $S_i, i=2,3,4, \dots,$ the symmetric, or permutation,  group of degree $i$.
	\item $D_i, i = 2,3,4,\dots,$ dihedral finite groups of degree $i$. (The symmetry groups of the planar regular $i$-gon.)
\end{itemize}

\subsection{The Representation of Finite Groups}
We will review some basic notations, concepts and results in the group
representation theory. Let $GL(n,\mathbb{C})$ be the group of $n
\times n$ (complex) non-singular matrices. A representation of $G$ of
degree $n$ is a homomorphism $\mathscr {D}$ from $G$ to
$GL(n,\mathbb{C})$.

For any two representations $\mathscr{D}_1$ and $\mathscr{D}_2$ of
$G$, respectively of degree $n_1$ and $n_2$, the direct sum
$\mathscr{D} = \mathscr{D}_1 \oplus \mathscr{D}_2$ can be defined as
another representation of degree $n_1+n_2$. For each $A \in G$,
$\mathscr{D}(A)$ is the $n_1+n_2$ invertible matrix 
$$ \mathscr{D}(A)=
\begin{pmatrix}
\mathscr{D}_1(A) & 0\\
0 & \mathscr{D}_2(A)\\
\end{pmatrix}.$$
Similarly, for representations
$\mathscr{D}_1,\mathscr{D}_2,\dots,\mathscr{D}_k$ of $G$ respectively
of degree $n_1, n_2, \dots,n_k,$ we can define a representation
$\mathscr{D} = \mathscr{D}_1 \oplus \mathscr{D}_2 \oplus \dots \oplus
\mathscr{D}_k$ of degree $n_1+n_2+\dots+n_k$. 

We say two representations $\mathscr{D}_1$ and $\mathscr{D}_2$ of
degree $n$ of $G$ are {\em equivalent}\/ if there is an invertible
$n \times n$ matrix $P$ such that for any $A \in G,$
$\mathscr{D}_2(A) = P\mathscr{D}_1(A)P^{-1}$.  A representation is
{\em reducible}\/ if it is equivalent to a direct sum of other
representations. Otherwise, we say that this representation is {\em
  irreducible}\/.

Let $\mathscr{D}$ be a representation of a group $G$ of degree $n$. A
character $\chi$ is a complex-valued function defined on group
$G$. For each $A \in G,$ let $\chi(A) = Tr(\mathscr{D}(A)).$ It is
easy to see that equivalent representations have the same
character. Most remarkably, representations which have the same
character are necessarily equivalent~\cite{2012Representation}.

The (Hermitian) inner product of characters can be defined as follows
$$(\chi_1,\chi_2) =\frac{1}{|G|} \sum_{A\in G} \overline{\chi_1}(A) \chi_2(A),$$ 
where $G$ is a finite group and $|G|$ is the number of elements in
$G$.  With this inner product, let $ \chi_1, \chi_2,\dots, \chi_k$
be characters of distinct (pair-wisely non-equivalent)
irreducible representations, we have, remarkably~\cite{2012Representation}
$$(\chi_i,\chi_j) = \delta_{ij}, \quad i, j = 1,2,\dots,k.$$

Suppose the finite group $G$ has $h$ conjugacy
  classes. Then the characters of representations of $G$ as vectors
span a $h$-dimensional complex vector space. To summarize, we state
the following theorem~\cite{2012Representation},  a classical result in
the representation theory of finite groups.

\begin{theorem} \label{Th1} For any finite group $G$ with $h$ conjugacy
  classes, there are $h$ non-equivalent irreducible
  representations. Let
  $\mathscr{D}_1, \mathscr{D}_2, \dots, \mathscr{D}_h$ be irreducible
  representations of group $G$. Let $\chi_1, \dots, \chi_h$ be the
  characters of these representations respectively. Then any
  representation $\mathscr{D}$ of $G$ with character $\chi$ is
  equivalent to
$$n_1 \mathscr{D}_1 \oplus n_2 \mathscr{D}_2 \oplus \dots \oplus n_h \mathscr{D}_h$$
where
$$n_i = (\chi,\chi_i)\in\mathbb{Z},~ i = 1, \dots ,h.$$
\end{theorem}

The above theorem provides an easy way to decompose any representation
$G$ into a direct sum of irreducible group representations. 

\subsection{The Dihedral Group $D_4$}
The dihedral group $D_4=<a,r~|~a^4=r^2=(ra)^2=e>$ has
$8$ elements. It is the symmetry group for a square. We
put $4$ identical particles at vertices and number them with
$1,~2,~3,~4$. Each element in $D_4$ acts as a permutation of the $4$
particles $\{1,~2,~3,~4\}$. 

\begin{table}[h]
\centering
\renewcommand{\tablename}{Table}
\renewcommand{\arraystretch}{1.2}  \doublerulesep 2.0pt
\begin{tabular}{|c|c|}
\hline
$D_4$ & Action\\
\hline
$e$ & Identity\\
\hline
$a$ & Rotating by $\frac{\pi}{2}$\\
\hline
$a^2$ & Rotating by $\pi$\\
\hline
$a^3$ & Rotating by $\frac{3\pi}{2}$\\
\hline
$r$ & Reflection\\
\hline
$ar$ & Reflection then Rotating by $\frac{\pi}{2}$\\
\hline
$a^2r$ & Reflection then rotating by $\pi$\\
\hline
$a^3r$ & Reflection then rotating by $\frac{3\pi}{2}$\\
\hline
\end{tabular}
\caption{The actions in $D_4$}
\label{action}
\end{table}

\begin{table}[h]
\centering
\renewcommand{\tablename}{Table}
\renewcommand{\arraystretch}{1.2}  \doublerulesep 2.0pt
\begin{tabular}{|c|c|c|c|c|c|}
\hline
$A$/ $\chi$ & $\chi_1$ & $\chi_2$ & $\chi_3$ & $\chi_4$ & $\chi_5$\\
\hline
$e$ & 1 & 1 & 1 & 1 & 2\\
\hline
$a,a^3$ & 1 & 1 & -1 & -1 & 0\\
\hline
$a^2$ & 1 & 1 & 1 & 1 & -2\\
\hline
$r,a^2$r & 1 & -1 & 1 & -1 & 0\\
\hline
$ar,a^3$r & 1 & -1 & -1 & 1 & 0\\
\hline 
\end{tabular}
\caption{\small The irreducible character table for $D_4$}
\label{d4}
\end{table}

The number of conjugacy classes for the group $D_4$ is $5$ which are
$\{e\},\{a,a^3\}$, $\{a^2\},\{r,a^2r\},\{ar,a^3r\}$. By Theorem
\ref{Th1}, the group $D_4$ also has $5$ irreducible
representations. The characters for the irreducible group
representation with $D_4$ are listed in Table \ref{d4}.  The degrees
of $\chi_1,\chi_2,\chi_3,\chi_4,\chi_5$ are $ 1,~1,~1,~1,~2$
respectively.

\subsection{The Dihedral group $D_3$, often known as, the Symmetric Group $S_3$}
The group $S_3$ has $6$ elements and is the symmetry group of
equilateral triangles. The elements of $S_3$ can be written as
matrices in the following 
$$ I=
\begin{pmatrix}
1 & 0 & 0\\
0 & 1 & 0\\
0 & 0 & 1\\
\end{pmatrix},
\quad
T=\begin{pmatrix}
1 & 0 & 0\\
0 & 0 & 1\\
0 & 1 & 0\\
\end{pmatrix},
\quad
R= \begin{pmatrix}
0 & 0 & 1\\
1 & 0 & 0\\
0 & 1 & 0\\
\end{pmatrix},   $$

$$ R^2=
\begin{pmatrix}
0 & 1 & 0\\
0 & 0 & 1\\
1 & 0 & 0\\
\end{pmatrix},
\quad
TR=\begin{pmatrix}
0 & 0 & 1\\
0 & 1 & 0\\
1 & 0 & 0\\
\end{pmatrix},
\quad
TR^2= \begin{pmatrix}
0 & 1 & 0\\
1 & 0 & 0\\
0 & 0 & 1\\
\end{pmatrix}.$$
This is actually a faithful representation of $S_3$ of degree $3.$ For $S_3$,
the conjugacy class are
$\left\{I\right\},\left\{R,R^2\right\},\left\{T,TR,TR^2\right\}$. So
there are three irreducible representations. First, for any element
$A$ in $S_3$, $\mathscr{D}_1(A) = 1$. This is the trivial
representation with degree $1$. Second, for any $A \in S_3$, let
$\mathscr{D}_2(A) = \mbox{det}(A)$. It is also a representation of
degree $1$. Third, there is a representation of degree $2$. It is 
$$ \mathscr{D}_3(I)=\bar{I}=
\begin{pmatrix}
1 & 0  \\
0 & 1  \\
\end{pmatrix},
\quad
\mathscr{D}_3(T)=\bar{T}=
\begin{pmatrix}
1 & 0  \\
0 & -1  \\
\end{pmatrix},$$

$$
\mathscr{D}_3(R)=\bar{R}=R( \frac{2\pi}{3})=
\begin{pmatrix}
-\frac{1}{2}       & -\frac{\sqrt{3}}{2}  \\
\frac{\sqrt{3}}{2} & -\frac{1}{2}  \\
\end{pmatrix},   $$

$$
\mathscr{D}_3(R^2)=\bar{R}^2=
\begin{pmatrix}
-\frac{1}{2} & \frac{\sqrt{3}}{2}  \\
-\frac{\sqrt{3}}{2} & -\frac{1}{2}  \\
\end{pmatrix},   $$

$$
\mathscr{D}_3(TR)=\bar{T}\bar{R}=
\begin{pmatrix}
-\frac{1}{2}       & -\frac{\sqrt{3}}{2}  \\
-\frac{\sqrt{3}}{2} & \frac{1}{2}  \\
\end{pmatrix}.   $$
It can be proved that $\mathscr{D}_1$, $\mathscr{D}_2$ and
$\mathscr{D}_3$ {enumerate} all the irreducible representations of
$S_3$. The character table for $S_3$ with
$\mathscr{D}_1$, $\mathscr{D}_2$ and $\mathscr{D}_3$ is presented in Table \ref{tab:my_label3}.

\begin{table}[h]
\centering
\renewcommand{\tablename}{Table}
\renewcommand{\arraystretch}{1.2}  \doublerulesep 2.0pt
\begin{tabular}{|c|c|c|c|}
\hline
A/ $\chi$ & $\chi_1$ & $\chi_2$ & $\chi_3$\\
\hline
$I$ & 1&1&2\\
\hline
$R, R^2$& 1&1&-1\\
\hline
$T,TR,TR^2$ & 1&-1&0\\
\hline %
\end{tabular}
\caption{The character table for $S_3$}
\label{tab:my_label3}
\end{table}

\section{APPLICATION}
In this section, we apply the representation theory of finite groups
to study central configurations with symmetries in the Newtonian
$n$-body problem. 

\subsection{The Technique from the Representation Theory}
Suppose $G$ is a finite group. Its element $A$ can act on some vector
space as a linear operator. Assume there is another linear operator
$H$ defined on the same vector space. We say that $H$ is {\em
  invariant under the group $G$}\/  if $HA=AH$ for each element $A$ in
$G$. Consider the eigenspace $V_\lambda$ for $H$ with $\phi \in
V_\lambda,$ we have 
$$H\phi=\lambda\phi \; \Rightarrow \; AH\phi=\lambda A\phi \; \Rightarrow \;
H(A\phi)=\lambda A\phi.$$ 
Then $A\phi$ is also an eigenvector for $H$ with the same
eigenvalue. Denote a basis for the eigenspace $V_\lambda$ by
$\phi_1,\phi_2,\dots,\phi_k$, it follows 
$$A\phi_j=\sum_{i=1}^{k}\mathscr{D}_{ij}(A)\phi_i, \quad j= 1,2,\dots,k.$$
Let $\mathscr{D}(A)$ be the matrix with entries $\mathscr{D}_{ij}(A).$
And if $A, B \in G,$ 
\begin{equation*}
\begin{aligned}
BA\phi_m &=\sum_{j=1}^{k}\mathscr{D}_{jm}(A)B\phi_j \\
&=\sum_{j=1}^{k}\mathscr{D}_{jm}(A)\sum_{i=1}^{k}\mathscr{D}_{ij}(B)\phi_i \\
&=\sum_{i=1}^{k}(\sum_{j=1}^{k}\mathscr{D}_{ij}(B)\mathscr{D}_{jm}(A))\phi_i, \quad m=1,2,\dots,k \\
\Rightarrow BA\phi_m &=\sum_{i=1}^{k}(\sum_{j=1}^{k}\mathscr{D}_{ij}(B)\mathscr{D}_{jm}(A))\phi_i \\
\Rightarrow \mathscr{D}(BA) &=\mathscr{D}(B)\mathscr{D}(A).
\end{aligned}
\end{equation*}

Therefore $\mathscr{D}$ is a group representation of $G$.

Suppose the $n \times n$ matrix $H$ is the Hessian for a smooth
function of $n$ variables. It is a symmetric matrix and has $n$
independent eigenvectors which form a basis of $\mathbb{R}^n$. In this
case, there is an invertible matrix $P$ (which acts as the change of
coordinates) such that $PAP^{-1}=\mathscr{D}(A)$ for all $A \in G$ and
$PHP^{-1}=H'$, where 
$$H'=
\begin{pmatrix}
\lambda_1 & 0 & \dots & 0\\
0 & \lambda_2 & \dots & 0\\
0 & 0 & \dots & 0\\
\vdots & \vdots & \ddots & \vdots\\
0 & 0 & \dots & \lambda_n\\
\end{pmatrix}.$$
Since the action of $A\in G$ does not intermix
{eigenspaces} of $H$. It is convenient to group the
eigenvectors with the same eigenvalue together and use them as a new
basis. The action of $G$, under the new basis, will be represented by a matrix of the form 
$$\mathscr{D}(A)=
\begin{pmatrix}
\mathscr{D}_1(A) & 0 & \dots & 0\\
0 & \mathscr{D}_2(A) & \dots & 0\\
0 & 0 & \dots & 0\\
\vdots & \vdots & \ddots & \vdots\\
0 & 0 & \dots & \mathscr{D}_k(A)\\
\end{pmatrix}$$ with each block represent an eigenspace of a single
eigenvalue and
$$\mbox{deg}(\mathscr{D}_1(A)) +\mbox{deg}(\mathscr{D}_2(A)) +
\dots+\mbox{deg} (\mathscr{D}_k(A))=n.$$
Then $H'\mathscr{D}(A)$ is simply

$$\begin{pmatrix}
\lambda_1\mathscr{D}_1(A) & 0 & \dots & 0\\
0 & \lambda_2\mathscr{D}_2(A) & \dots & 0\\
0 & 0 & \dots & 0\\
\vdots & \vdots & \ddots & \vdots\\
0 & 0 & \dots & \lambda_k\mathscr{D}_k(A)\\
\end{pmatrix}.$$

This implies that $\mathscr{D}$ is equivalent to
$\mathscr{D}_1 \oplus \mathscr{D}_2 \oplus \dots \oplus
\mathscr{D}_k.$ In general, the constituents $\mathscr{D}_i$,
{$i=1,2,\dots,k$, associated with each eigenspace might
  not be irreducible. However, it can be further decomposed into
  irreducible constituents.  This can be easily done by choosing proper
  coordinates in a given eigenspace. We remark that, in a typical
  case, these constituents are already irreducible. However, there
  may be accidental degeneracies (perhaps due to some other
  symmetries) such that two sets of eigenvectors have the same
  eigenvalue. Indeed, in the $n$-body example we are considering, we
  do have a symmetry that forces two zero eigenvalues that correspond
  to two different irreducible representations. In any case, the degeneracy is not
  an issue here, since we already know the irreducible constituents of
  $\mathscr{D}$ by the character computation. We may assume all
  $\mathscr{D}_i$, $i =1, \ldots, k$ are irreducible, even though
  their corresponding eigenvalues may not necessarily be distinct.}

  To summarize, let $G$ be a finite group, we can find all irreducible
  representations of $G$ and their characters. This allows us to
  decompose any group representation $\mathscr{D}$ of $G$ into
  irreducible components. If a symmetric matrix $H$ is invaraint under action
  $\mathscr{D}$, this decomposition of $\mathscr{D}$ can be done by properly choosing
  eigenvectors as new basis.
  We can
  assume that, under the new coordinates, $\mathscr{D}_i$, for each $i$, is an irreducible group
  representation $\mathscr{D}$. Finally, for each element $A$ in $G$, we can compute
  the trace of $\mathscr{D}_i(A), ~i=1,2,\dots,k$. Since
  $Tr(H\mathscr{D}(A))=Tr(PH\mathscr{D}(A)P^{-1})=Tr(H'\mathscr{D}(A))$,
  we obtain a list of simple linear equations involving the
  eigenvalues of $H$. These equations may or may not be sufficient in
  solving the eigenvalue problem for $H$, they certainly help. Indeed,
  with some other information on the matrix $H$, an otherwise
  impossible problem may become algebraically solvable.

  \medskip

  In the next subsections, we use the above technique to study two
  examples in the four body problems.

\subsection{The Square Configuration}
\begin{figure}[ht]
	\centering
	\includegraphics[width=0.7\linewidth]{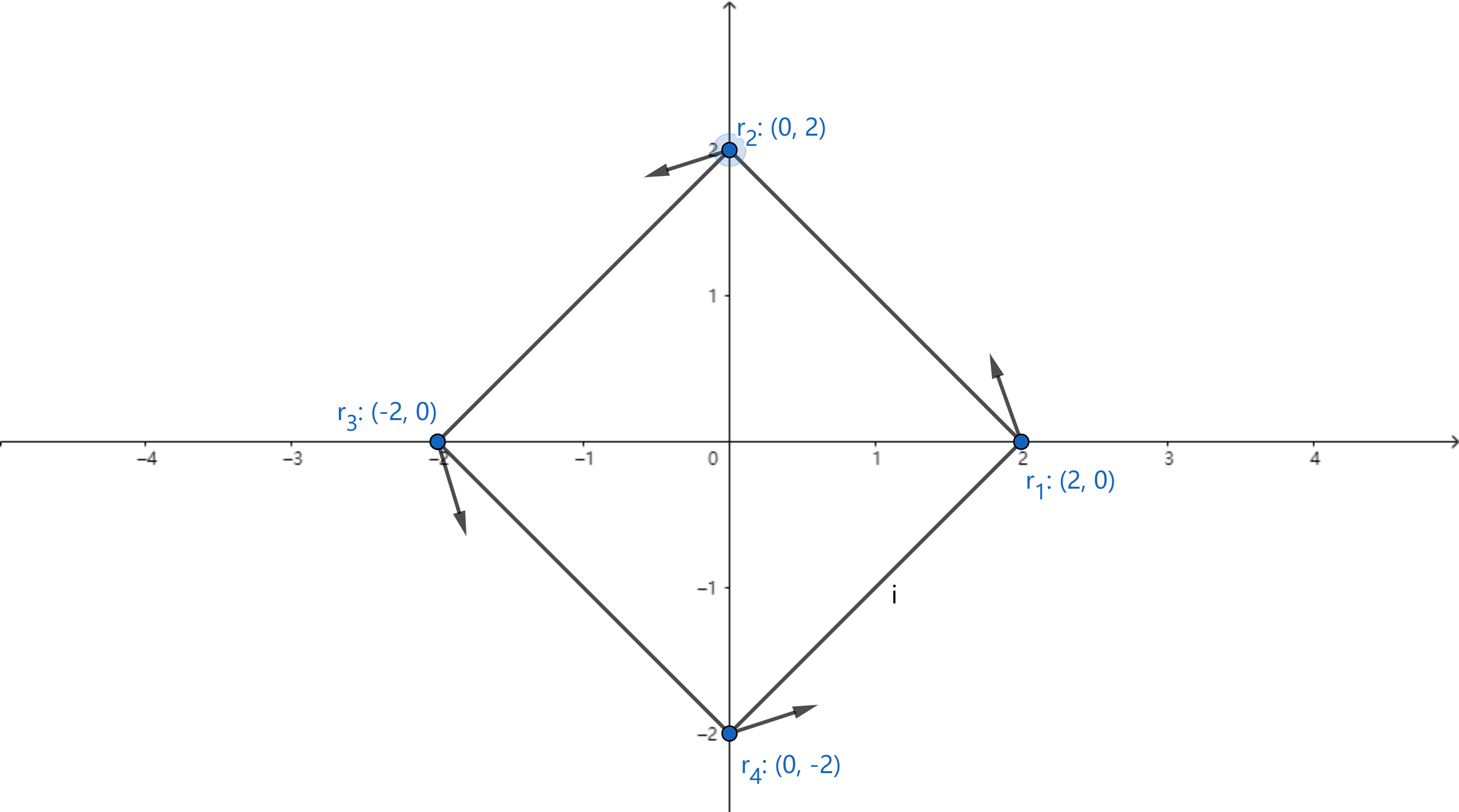}
	\renewcommand{\figurename}{Figure}
	\caption{The Square Configuration}
	\label{fig:figure1}
\end{figure}

The first example is the $4$-body problem in $\mathbb{R}^2$ with four equal masses.
The central configurations are critical points of the function
$\sqrt{I}U$, where 
$$ I=\sum_{i=1}^{4} \frac{1}{2} (x_i^2+y_i^2),$$
$$U=\sum_{1\le i<j \le 4} \frac{1}{\sqrt{(x_i-x_j)^2+(y_i-y_j)^2}}. $$
Then square centered at the origin, as shown in Figure \ref{fig:figure1}, is a central configuration.

Let $z \in \mathbb{R}^8$ be the vector
$z=(x_1,y_1,x_2,y_2,x_3,y_3,x_4,y_4).$ Then the configuration for the
square is 
$$z_0=(2,0,0,2,-2,0,0,-2).$$
Clearly $z_0$ is a critical point for $\sqrt{I}U,$ and the Hessian
$ \frac{\partial^2(\sqrt{2I}U)}{\partial z^2}$ at $z_0$ is
$$H_1=
\begin{pmatrix}
\frac{3\sqrt{2}}{16}+\frac{5}{32} & 0 & -\frac{\sqrt{2}}{16} & -\frac{3}{32} & \frac{3\sqrt{2}}{16}-\frac{1}{32} & 0 & -\frac{\sqrt{2}}{16} & \frac{3}{32} \\
0 & \frac{3\sqrt{2}}{8}+\frac{1}{16} & \frac{3\sqrt{2}}{16} & -\frac{\sqrt{2}}{16} & 0 & \frac{1}{16} & -\frac{3\sqrt{2}}{16} & -\frac{\sqrt{2}}{16}\\
-\frac{\sqrt{2}}{16} & \frac{3\sqrt{2}}{16} & \frac{3\sqrt{2}}{8}+\frac{1}{16} & 0 & -\frac{\sqrt{2}}{16} & -\frac{3\sqrt{2}}{16} & \frac{1}{16} & 0 \\
-\frac{3}{32} & -\frac{\sqrt{2}}{16} & 0 & \frac{3\sqrt{2}}{16}+\frac{5}{32} & \frac{3}{32} & -\frac{\sqrt{2}}{16} & 0 & \frac{3\sqrt{2}}{16}-\frac{1}{32}\\
\frac{3\sqrt{2}}{16}-\frac{1}{32} & 0 & -\frac{\sqrt{2}}{16} & \frac{3}{32} & \frac{3\sqrt{2}}{16}+\frac{5}{32} & 0 & -\frac{\sqrt{2}}{16} & -\frac{3}{32}\\
0 & \frac{1}{16} & -\frac{3\sqrt{2}}{16} & -\frac{\sqrt{2}}{16} & 0 & \frac{3\sqrt{2}}{8}+\frac{1}{16} & \frac{3\sqrt{2}}{16} & -\frac{\sqrt{2}}{16}\\
-\frac{\sqrt{2}}{16} & -\frac{3\sqrt{2}}{16} & \frac{1}{16} & 0 & -\frac{\sqrt{2}}{16} & \frac{3\sqrt{2}}{16} & \frac{3\sqrt{2}}{8}+\frac{1}{16} & 0\\
\frac{3}{32} & -\frac{\sqrt{2}}{16} & 0 & \frac{3\sqrt{2}}{16}-\frac{1}{32} & -\frac{3}{32} & -\frac{\sqrt{2}}{16} & 0 & \frac{3\sqrt{2}}{16}+\frac{5}{32}
\end{pmatrix}.$$
\par
We now use the dihedral group $D_4$ to find the eigenvalues of $H_1.$
As shown in Figure \ref{fig:figure1}, the square has four equal
{masses} at its vertices. Every element in $D_4$
{acts} on the square by permuting the particles.

\begin{figure}[ht]
	\centering
	\includegraphics[width=0.7\linewidth]{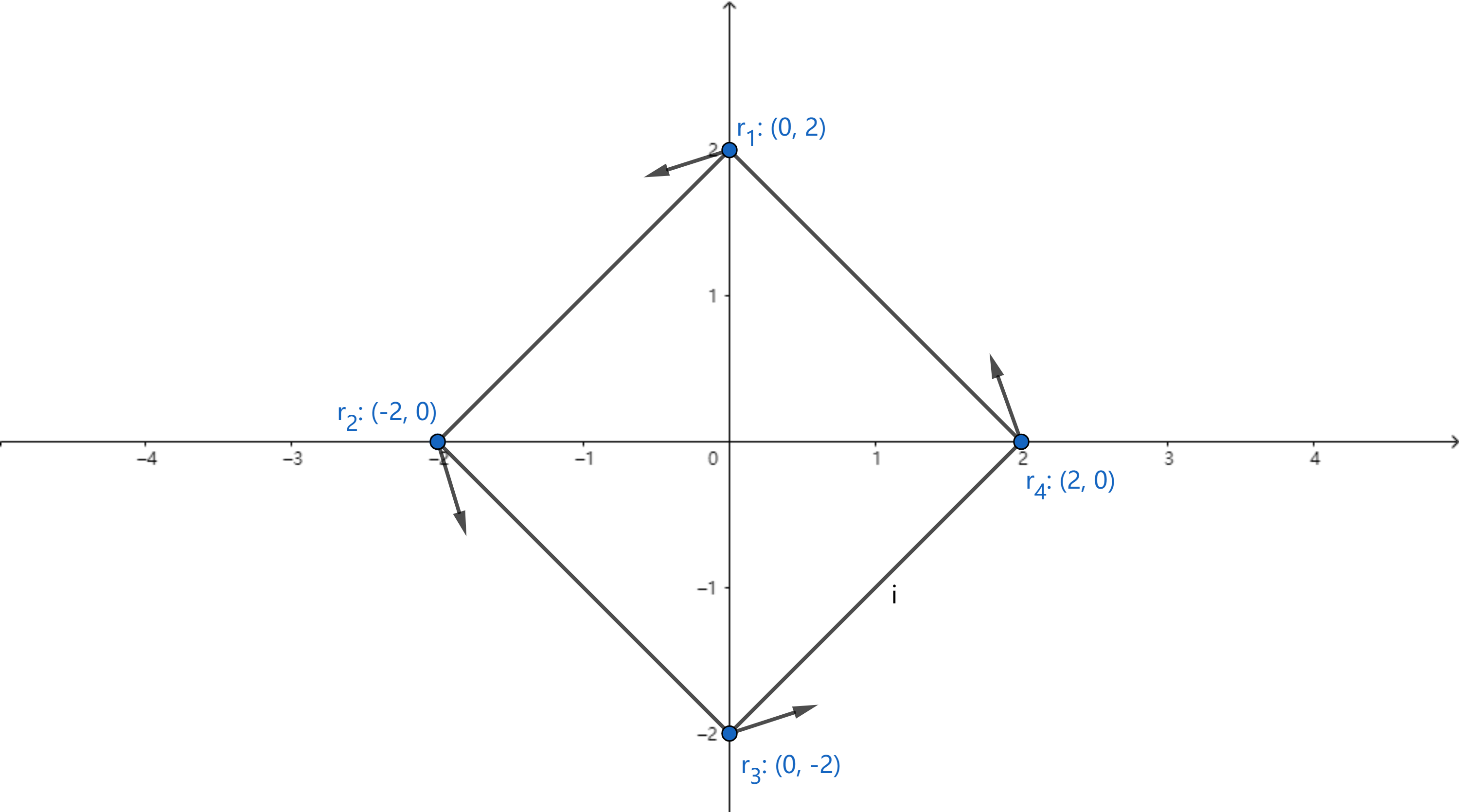}
	\renewcommand{\figurename}{Figure}
	\caption{The square after rotating by $\frac{\pi}{2}$}
	\label{fig:figure2}
\end{figure}

For example, take $a\in D_4$ that rotates square by
$\frac{\pi}{2}$, the position of individual masses become 
$$\begin{pmatrix}
x_1'\\y_1'
\end{pmatrix}
=R(\frac{\pi}{2})
\begin{pmatrix}
 x_4\\y_4
\end{pmatrix}, \quad
\begin{pmatrix}
x_2'\\y_2'
\end{pmatrix}
=R(\frac{\pi}{2})
\begin{pmatrix}
x_1\\y_1
\end{pmatrix},$$
$$\begin{pmatrix}
x_3'\\y_3'
\end{pmatrix}
=R(\frac{\pi}{2})
\begin{pmatrix}
x_2\\y_2
\end{pmatrix}, \quad
\begin{pmatrix}
x_4'\\y_4'
\end{pmatrix}
=R(\frac{\pi}{2})
\begin{pmatrix}
x_3\\y_3
\end{pmatrix}.$$
The four body central configuration shown in Figure \ref{fig:figure1}
{transforms to}  Figure \ref{fig:figure2}.

With the vector $z=(x_1,y_1,x_2,y_2,x_3,y_3,x_4,y_4)$, we have
$$z'=
\begin{pmatrix}
0 & 0 & 0 & R(\frac{\pi}{2}) \\
R(\frac{\pi}{2}) & 0 & 0 & 0 \\
0 & R(\frac{\pi}{2}) & 0 & 0 \\
0 & 0 & R(\frac{\pi}{2}) & 0 \\
\end{pmatrix}z,$$
where $z'$ is the vector after the transformation. Hence, action $a$ can be represented by the $8 \times 8$ matrix

$$\mathscr{D}(a)=
\begin{pmatrix}
0 & 0 & 0 & R(\frac{\pi}{2}) \\
R(\frac{\pi}{2}) & 0 & 0 & 0 \\
0 & R(\frac{\pi}{2}) & 0 & 0 \\
0 & 0 & R(\frac{\pi}{2}) & 0 \\
\end{pmatrix}, \quad \mbox{where} \quad R(\frac{\pi}{2})=
\begin{pmatrix}
0 & -1\\
1 & 0
\end{pmatrix}. $$
Similarly $a^2, r, ar \in D_4$ are respectively rotating $\pi$, reflection and rotating $\frac{\pi}{2}$ after reflection. They have representations as follows
$$\mathscr{D}(a^2)=
\begin{pmatrix}
0 & 0 & R(\pi) & 0 \\
0 & 0 & 0 & R(\pi) \\
R(\pi) & 0 & 0 & 0 \\
0 & R(\pi) & 0 & 0 \\
\end{pmatrix}, \quad \mbox{where} \quad R(\pi)=
\begin{pmatrix}
-1 & 0\\
0 & -1
\end{pmatrix}, $$
$$\mathscr{D}(r)=
\begin{pmatrix}
0 & 0 & F & 0 \\
0 & F & 0 & 0 \\
F & 0 & 0 & 0 \\
0 & 0 & 0 & F \\
\end{pmatrix}, \quad \mbox{where} \quad F=
\begin{pmatrix}
-1 & 0\\
0 & 1
\end{pmatrix}, $$
$$\mathscr{D}(ar)=
\begin{pmatrix}
0 & 0 & 0 & G \\
0 & 0 & G & 0 \\
0 & G & 0 & 0 \\
G & 0 & 0 & 0 \\
\end{pmatrix}, \quad \mbox{where} \quad G=R(\frac{\pi}{2})F=
\begin{pmatrix}
0 & -1\\
-1 & 0
\end{pmatrix}. $$
\par
Since the character of $\mathscr{D}$ is constant on conjugacy classes,
its character can be easily calculated. 
$$\chi(e)=8, \quad \chi(a)=0,$$
$$\chi(a^2)=0, \quad \chi(r)=0,\quad \chi(ar)=0.$$

\begin{table}[h]
\centering
\renewcommand{\tablename}{Table}
\renewcommand{\arraystretch}{1.2}  \doublerulesep 2.0pt
\begin{tabular}{|c|c|c|c|c|c|c|}
\hline
A & $\chi$ & $\chi_1$ & $\chi_2$ & $\chi_3$ & $\chi_4$ & $\chi_5$\\
\hline
e & 8 & 1 & 1 & 1 & 1 & 2\\
\hline
a,a$^3$ & 0 & 1 & 1 & -1 & -1 & 0\\
\hline
a$^2$ & 0 & 1 & 1 & 1 & 1 & -2\\
\hline
r,a$^2$r & 0 & 1 & -1 & 1 & -1 & 0\\
\hline
ar,a$^3$r & 0 & 1 & -1 & -1 & 1 & 0\\
\hline
\end{tabular}
\caption{The character of $\mathscr{D}$ and irreducible characters for $D_4$}
\label{tab:my_label4}
\end{table}

Now, we can decompose $\mathscr{D}$ into irreducible
representations.  According to Theorem \ref{Th1}, we have
$$
\begin{aligned}
n_1=(\chi,\chi_1)&=\frac{1}{8}\times 8=1, \quad n_2=(\chi,\chi_2)=1, \\
n_3=(\chi,\chi_3)&=1,\quad n_4=(\chi,\chi_4)=1,\quad n_5=(\chi,\chi_5)=2.
\end{aligned}
$$
Hence $\chi=\chi_1+\chi_2+\chi_3+\chi_4+2\chi_5,$ we conclude
that  $$\mathscr{D} \sim
\mathscr{D}_1 \oplus \mathscr{D}_2 \oplus \mathscr{D}_3
\oplus \mathscr{D}_4 \oplus 2\mathscr{D}_5.$$ 

On the other hand, the Hessian $H_1$ is invariant under the action of the symmetry group $D_4$,
i.e., $\mathscr{D}(A)H_1=H_1\mathscr{D}(A)$ for every element $A \in D_4.$
This implies that the action $\mathscr{D}(A)$ for any $A \in D_4$
doesn't intermix the eigenspaces of $H_1$. As discussed in the last
subsection, with a proper choice of eigenvectors of $H_1$ as our new
coordinate system, the representation of
$\mathscr{D}$ of dihedral group $D_4$ is now given by 

$$\mathscr{D}'(A)=
\begin{pmatrix}
\mathscr{D}_1(A) & 0 & 0 & 0 & 0 & 0 \\
0 & \mathscr{D}_2(A) & 0 & 0 & 0 & 0 \\
0 & 0 & \mathscr{D}_3(A) & 0 & 0 & 0 \\
0 & 0 & 0 & \mathscr{D}_4(A) & 0 & 0 \\
0 & 0 & 0 & 0 & \mathscr{D}_5(A) & 0 \\
0 & 0 & 0 & 0 & 0 & \mathscr{D}_5(A) \\
\end{pmatrix}. $$
{The degrees of
  $\mathscr{D}_1,\mathscr{D}_2,\mathscr{D}_3,\mathscr{D}_4,\mathscr{D}_5$
  are, respectively, $1,~1,~1,~1,~2,$}, therefore the matrix $H_1'$
transformed from $H_1$ has the form
$$H_1'=
\begin{pmatrix}
\lambda_1 & 0 & 0 & 0 & 0 & 0 & 0 & 0\\
0 & \lambda_2 & 0 & 0 & 0 & 0 & 0 & 0\\
0 & 0 & \lambda_3 & 0 & 0 & 0 & 0 & 0\\
0 & 0 & 0 & \lambda_4 & 0 & 0 & 0 & 0\\
0 & 0 & 0 & 0 & \lambda_5 & 0 & 0 & 0\\
0 & 0 & 0 & 0 & 0 & \lambda_5 & 0 & 0\\
0 & 0 & 0 & 0 & 0 & 0 & \lambda_6 & 0\\
0 & 0 & 0 & 0 & 0 & 0 & 0 & \lambda_6\\
\end{pmatrix}. $$
The trace of $H_1\mathscr{D}(A)$ is therefore equivalent to the trace
of 
$$\begin{pmatrix}
\lambda_1\mathscr{D}_1(A) & 0 & 0 & 0 & 0 & 0 \\
0 & \lambda_2\mathscr{D}_2(A) & 0 & 0 & 0 & 0 \\
0 & 0 & \lambda_3\mathscr{D}_3(A) & 0 & 0 & 0 \\
0 & 0 & 0 & \lambda_4\mathscr{D}_4(A) & 0 & 0 \\
0 & 0 & 0 & 0 & \lambda_5\mathscr{D}_5(A) & 0 \\
0 & 0 & 0 & 0 & 0 & \lambda_6\mathscr{D}_5(A) \\
\end{pmatrix}, $$

i.e.,
\begin{equation}
	Tr(H_1\mathscr{D}(A))=\lambda_1\chi_1(A)+\lambda_2\chi_2(A) +
	\lambda_3\chi_3(A) + \lambda_4\chi_4(A) + \lambda_5\chi_5(A) +
	\lambda_6\chi_5(A). \label{trace d4}
\end{equation}
On the other hand, the trace of $H_1\mathscr{D}(A)$ can be calculated
in the original coordinate. Then we have the following equations
\begin{eqnarray}
Tr(H_1\mathscr{D}(e)) &=& Tr(H_1) =
                          \frac{9\sqrt{2}}{4}+\frac{7}{8} = \lambda_1+
                          \lambda_2+\lambda_3+\lambda_4+ 2(\lambda_5+ \lambda_6),
  \nonumber \\
Tr(H_1\mathscr{D}(a)) &=&
                          -\frac{3}{8}-\frac{3\sqrt{2}}{4}=\lambda_1+\lambda_2-
                          \lambda_3-\lambda_4, \nonumber\\
Tr(H_1\mathscr{D}(a^2)) &=&  -\frac{1}{8}- \frac{3\sqrt{2}}{4} =
                            \lambda_1 +\lambda_2+\lambda_3+
                            \lambda_4-2(\lambda_5+\lambda_6),
                            \nonumber \\
Tr(H_1\mathscr{D}(r)) &=&
                          \frac{3}{8}-\frac{3\sqrt{2}}{4}=\lambda_1-\lambda_2+\lambda_3-\lambda_4,
                          \nonumber \\
Tr(H_1\mathscr{D}(ar)) &=&
                           -\frac{3}{8}+\frac{3\sqrt{2}}{4}=\lambda_1-\lambda_2-\lambda_3+\lambda_4. \nonumber
\end{eqnarray}

There are 5 equations for 6 unknowns. Solving these five equations, we obtain
$$\lambda_1=\lambda_2=0, \quad \lambda_3=\frac{3}{8}, \quad
\lambda_4=\frac{3\sqrt{2}}{4}, \quad
\lambda_5+\lambda_6=\frac{3\sqrt{2}}{4}+\frac{1}{4}.$$ 
We need one more equation to solve for all 6 eigenvalues. This is
easy.  Observe that the determinant of $H+I_8$ can be easily
computed, we have $$\frac{11}{8}(\frac{3\sqrt{2}}{4}+1)(\lambda_5+1)^2(\lambda_6+1)^2=\mbox{det}(H+I_8)=\frac{340505}{32768} + \frac{963897\sqrt{2}}{131072}$$
$$\Rightarrow	(\lambda_5+1)(\lambda_6+1)=\frac{97}{64}+ \frac{27\sqrt{2}}{32}, $$
we get
$$\lambda_1=\lambda_2=0, \quad \lambda_3=\frac{3}{8}, \quad \lambda_4=\frac{3\sqrt{2}}{4}, $$
$$\lambda_5=\frac{\sqrt{2}}{4}+\frac{1}{8}, \quad \lambda_6=\frac{\sqrt{2}}{2}+\frac{1}{8}.$$

We point out that two zero eigenvalues are expected, as they correspond
to rotation and scaling invariance of  $\sqrt{I} U$.

\subsection{The Center + Equilateral Triangle Configuration}
\begin{figure}[ht]
	\centering
	\includegraphics[width=0.7\linewidth]{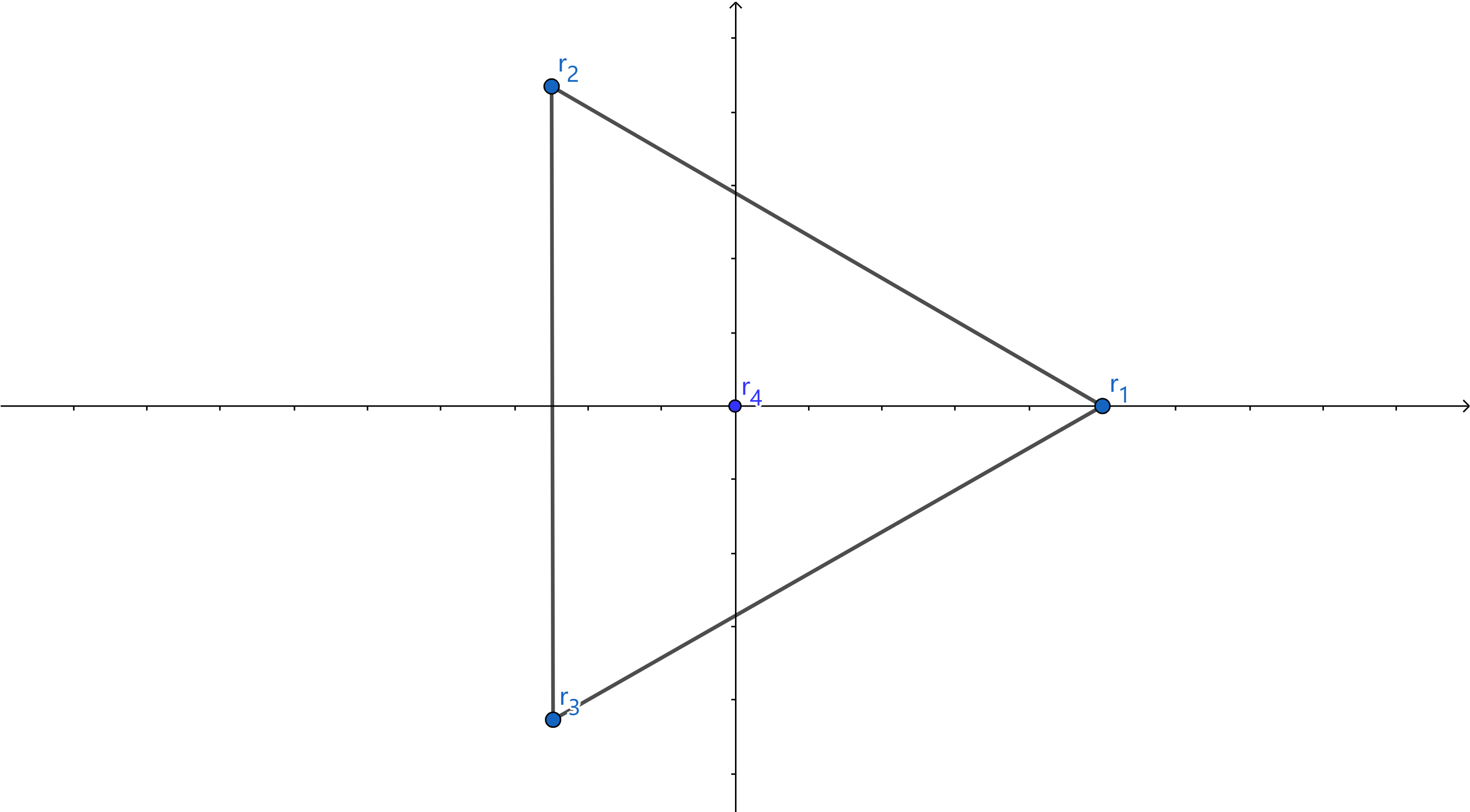}
	\renewcommand{\figurename}{Figure}
	\caption{The center + equilateral triangle solution for the $4$ body problem}
	\label{fig:figure3}
\end{figure}

Next we consider the central configuration with three {particles} with
mass 1 at the vertices of an {equilateral} triangle and a fourth particle
with mass $m$ at the origin, as shown in Figure \ref{fig:figure3}.
Clearly, this is a central configuration. 

Let 
$$ I=\frac{1}{2}(\sum_{i=1}^{3}(x_i^2+y_i^2)+ m(x_4^2+y_4^2)),$$
$$U=\sum_{1\le i<j \le 3}
\frac{1}{\sqrt{(x_i-x_j)^2+(y_i-y_j)^2}}+m\sum_{i=1}^{3}
\frac{1}{\sqrt{(x_i-x_4)^2+(y_i-y_4)^2}}. $$ 
Let $z \in \mathbb{R}^8$ be the vector
$z=(x_1,y_1,x_2,y_2,x_3,y_3,x_4,y_4).$ Then the equilateral {triangle}
can be described by 
$$z_0=(1,0,-\frac{1}{2},\frac{\sqrt{3}}{2},-\frac{1}{2},- \frac{\sqrt{3}}{2},0,0)$$
This is a critical point of the function $\sqrt{I}U$. The Hessian $H_2$
for $ \frac{\partial^2(\sqrt{2I}U)}{\partial z^2}$ at $z_0$ is: 
$${\footnotesize
\begin{pmatrix}
\frac{5}{6}+2\sqrt{3}m & 0 & \frac{1}{12}+\frac{\sqrt{3}m}{2} & -\frac{\sqrt{3}}{4}-\frac{3m}{2} & \frac{1}{12}+\frac{\sqrt{3}m}{2} & \frac{\sqrt{3}}{4}+\frac{3m}{2} & -2\sqrt{3}m & 0 \\
0 & \frac{5}{6} & \frac{\sqrt{3}}{4} & \frac{1}{12} & -\frac{\sqrt{3}}{4} & \frac{1}{12} & 0 & \sqrt{3}m \\
\frac{1}{12}+\frac{\sqrt{3}m}{2} & \frac{\sqrt{3}}{4} & \frac{5}{6}+\frac{\sqrt{3}m}{2} & -\frac{3m}{2} & \frac{1}{12}-\frac{\sqrt{3}m}{4} & -\frac{\sqrt{3}}{4}-\frac{3m}{4} & \frac{\sqrt{3}m}{4} & \frac{9m}{4} \\
-\frac{\sqrt{3}}{4}-\frac{3m}{2} & \frac{1}{12} & -\frac{3m}{2} & \frac{5}{6}+\frac{3\sqrt{3}m}{2} & \frac{\sqrt{3}}{4}+\frac{3m}{4} & \frac{1}{12}+\frac{3\sqrt{3}m}{4} & \frac{9m}{4} & -\frac{5\sqrt{3}m}{4} \\
\frac{1}{12}+\frac{\sqrt{3}m}{2} & -\frac{\sqrt{3}}{4} & \frac{1}{12}-\frac{\sqrt{3}m}{4} & \frac{\sqrt{3}}{4}+\frac{3m}{4} & \frac{5}{6}+\frac{\sqrt{3}m}{2} & \frac{3m}{2} & \frac{\sqrt{3}m}{4} & -\frac{9m}{4}\\
\frac{\sqrt{3}}{4}+\frac{3m}{2} & \frac{1}{12} & -\frac{\sqrt{3}}{4}-\frac{3m}{4} & \frac{1}{12}+\frac{3\sqrt{3}m}{4} & \frac{3m}{2} & \frac{5}{6}+\frac{3\sqrt{3}m}{2}  & -\frac{9m}{4} & -\frac{5\sqrt{3}}{4}\\
-2\sqrt{3}m & 0 & \frac{\sqrt{3}m}{4} & \frac{9m}{4} & \frac{\sqrt{3}m}{4} & -\frac{9m}{4} & \frac{\sqrt{3}(2m^2+3m)}{2}+m & 0\\
0 & \sqrt{3}m & \frac{9m}{4} & -\frac{5\sqrt{3}m}{4} & -\frac{9m}{4} & -\frac{5\sqrt{3}m}{4} & 0 & \frac{\sqrt{3}(2m^2+3m)}{2}+m \\
\end{pmatrix}}$$
The group $S_3$ is a symmetry group for the central configuration. An
action $A$ in $S_3$ displaces particles by a linear transformations
$\mathscr{D}(A)$. For instance, the action $R$ that rotates the entire
configuration by $\frac{2\pi}{3},$ {can} be described by
$$\begin{pmatrix}
x_1'\\y_1'
\end{pmatrix}
=R(\frac{2\pi}{3})
\begin{pmatrix}
 x_3\\y_3
\end{pmatrix}, \quad
\begin{pmatrix}
x_2'\\y_2'
\end{pmatrix}
=R(\frac{2\pi}{3})
\begin{pmatrix}
x_1\\y_1
\end{pmatrix},$$
$$\begin{pmatrix}
x_3'\\y_3'
\end{pmatrix}
=R(\frac{2\pi}{3})
\begin{pmatrix}
x_2\\y_2
\end{pmatrix}, \quad
\begin{pmatrix}
x_4'\\y_4'
\end{pmatrix}
=R(\frac{2\pi}{3})
\begin{pmatrix}
x_4\\y_4
\end{pmatrix}.$$
It induces a representation for $R$,
$$\mathscr{D}(R)=
\begin{pmatrix}
	0 & 0 & \bar{R} & 0  \\
	\bar{R} & 0 & 0 & 0  \\
	0 & \bar{R} & 0 & 0 \\
	0 & 0 & 0 & \bar{R} \\
\end{pmatrix},
\quad \mbox{where} \quad
\bar{R} = R(\frac{2\pi}{3})=
\begin{pmatrix}
	-\frac{1}{2}  & -\frac{\sqrt{3}}{2} \\
	\frac{\sqrt{3}}{2} & -\frac{1}{2} \\
\end{pmatrix}. $$
Similarly, we have
$$
\mathscr{D}(T)=
\begin{pmatrix}
\bar{T} & 0 & 0 & 0 \\
0 & 0 & \bar{T} & 0  \\
0 & \bar{T} & 0 & 0  \\
0 & 0 & 0 & \bar{T} \\
\end{pmatrix},
\quad \mbox{where} \quad
\bar{T} =
\begin{pmatrix}
1  & 0 \\
0 & -1 \\
\end{pmatrix}. $$

Therefore, $\mathscr{D}$ is a representation of $S_3$ with degree $8.$ The character of $\mathscr{D}$ can be calculated with
$$\chi(I)=8, \quad \chi(R)=-1, \quad \chi(T)=0.$$
Referring the character table for $S_3,$ we have Table \ref{tab:my_label5}.
\begin{table}[ht]
\centering
\renewcommand{\tablename}{Table}
\renewcommand{\arraystretch}{1.2}  \doublerulesep 2.0pt
\begin{tabular}{|c|c|c|c|c|}
	\hline
	$A$ & $\chi$ & $\chi_1$ & $\chi_2$ & $\chi_3$\\
	\hline
	$I$ & 8 & 1 & 1 & 2\\
	\hline
	$R$, $R^2$& -1 & 1&1&-1\\
	\hline
	$T$,$TR$,$TR^2$ & 0 & 1&-1&0\\
	\hline 
\end{tabular}
\caption{\small The character of $\mathscr{D}$ and character table for $S_3$}
    \label{tab:my_label5}
\end{table}

According to Theorem \ref{Th1}, it implies that
$$
\begin{aligned}
n_1=(\chi,\chi_1) &=\frac{1}{6}\times (8 \times 1-2\times 1\times 1 + 3\times 0\times 1) = 1,   \\
n_2=(\chi,\chi_2) &=1, \quad n_3=(\chi,\chi_3)=3.
\end{aligned}
$$
This implies that $\chi=\chi_1+\chi_2+3\chi_3$. Therefore
$\mathscr{D}$ is equivalent to $\mathscr{D}_1 \oplus \mathscr{D}_2
\oplus 3\mathscr{D}_3.$  It can be verified that
$H_2\mathscr{D}(A)=\mathscr{D}(A)H_2$ for all $A \in S_3.$ Since $H_2$
is invariant under the symmetry group $S_3$ , the action by
$\mathscr{D}(A)$ does not {intermix the eigenspaces of
  $H_2$.} Ahain, by choosing proper eigenvectors for $H_2$ as a basis in $\mathbb{R}^8$, there is an invertible matrix $P$ such that $P\mathscr{D}(A)P^{-1}=\mathscr{D}(A)',$ where $\mathscr{D}(A)'$ is equivalent to
$$\begin{pmatrix}
\mathscr{D}_1(A) & 0 & 0 & 0 & 0  \\
0 & \mathscr{D}_2(A) & 0 & 0 & 0  \\
0 & 0 & \mathscr{D}_3(A) & 0 & 0  \\
0 & 0 & 0 & \mathscr{D}_3(A) & 0  \\
0 & 0 & 0 & 0 & \mathscr{D}_3(A)  \\
\end{pmatrix}, $$
and degrees of $\mathscr{D}_1, \mathscr{D}_2,\mathscr{D}_3$ are
respectively $1,~1,~2$. Under the new coordinates,  $H_2$ takes the form
$$H_2'=
\begin{pmatrix}
\lambda_1 & 0 & 0 & 0 & 0 & 0 & 0 & 0\\
0 & \lambda_2 & 0 & 0 & 0 & 0 & 0 & 0\\
0 & 0 & \lambda_3 & 0 & 0 & 0 & 0 & 0\\
0 & 0 & 0 & \lambda_3 & 0 & 0 & 0 & 0\\
0 & 0 & 0 & 0 & \lambda_4 & 0 & 0 & 0\\
0 & 0 & 0 & 0 & 0 & \lambda_4 & 0 & 0\\
0 & 0 & 0 & 0 & 0 & 0 & \lambda_5 & 0\\
0 & 0 & 0 & 0 & 0 & 0 & 0 & \lambda_5\\
\end{pmatrix}. $$
Therefore the trace of $H_2\mathscr{D}(A)$ is equal to the trace of
$$\begin{pmatrix}
\lambda_1\mathscr{D}_1(A) & 0 & 0 & 0 & 0  \\
0 & \lambda_2\mathscr{D}_2(A) & 0 & 0 & 0  \\
0 & 0 & \lambda_3\mathscr{D}_3(A) & 0 & 0  \\
0 & 0 & 0 & \lambda_4\mathscr{D}_3(A) & 0  \\
0 & 0 & 0 & 0 & \lambda_5\mathscr{D}_3(A)  \\
\end{pmatrix}. $$
We obtain a linear equations on eigenvalues of $H_2$:
\begin{equation}
Tr(H_2\mathscr{D}(A))=\lambda_1\chi_1(A)+\lambda_2\chi_2(A)+\lambda_3\chi_3(A)+\lambda_4\chi_3(A)+\lambda_5\chi_3(A). \label{character}
\end{equation}
For $A=I,$
$$\lambda_1+\lambda_2+2(\lambda_3+\lambda_4+\lambda_5)=Tr(H_2)=2\sqrt{3}m^2+(9\sqrt{3}+2)m+5.$$
For $A=R,$
$$\lambda_1+\lambda_2-(\lambda_3+\lambda_4+\lambda_5)=Tr(H_2\mathscr{D}(R))=-\frac{1}{2}(2\sqrt{3}m^2+(9\sqrt{3}+2)m+5).$$
For $A=T,$
$$\lambda_1-\lambda_2+0(\lambda_3+\lambda_4+\lambda_5)=Tr(H_2\mathscr{D}(T))=0.$$
Solving these equations, we have
$$\lambda_1=\lambda_2=0, \quad \lambda_3+\lambda_4+\lambda_5=\frac{1}{2}(2\sqrt{3}m^2+(9\sqrt{3}+2)m+5).$$
Again, we don't have enough equations from the symmetry to find all
eigenvalues. We need to find some additional equations. Now consider the characteristic polynomial
$$f(\lambda)=|\lambda I_8-H|=\lambda^8+a_1\lambda^7+\dots+a_8.$$
Since $\lambda_1=\lambda_2=0,$ we have
$$\lambda_3^2+\lambda_4^2+\lambda_5^2+4\lambda_3\lambda_4+4\lambda_3\lambda_5+4\lambda_4\lambda_5=a_2,$$
$$\lambda_3^2\lambda_4^2\lambda_5^2=a_6.$$
Then
$$
\lambda_3\lambda_4+\lambda_3\lambda_5+\lambda_4\lambda_5 =\frac{1}{2}(a_2-(\lambda_3+\lambda_4+\lambda_5)^2), \quad
\lambda_3\lambda_4\lambda_5=\sqrt{a_6}. $$

\begin{figure}[h]
	\centering
	\includegraphics[width=0.5\linewidth]{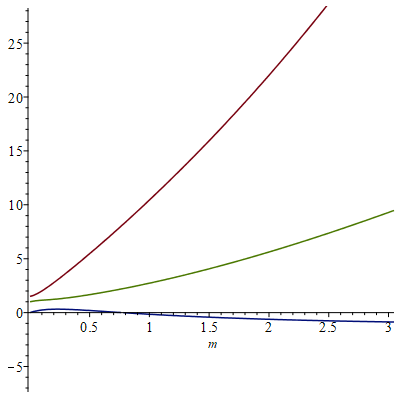}
	\renewcommand{\figurename}{Figure}
	\caption{Eigenvalues with the parameter $m$}
	\label{fig:figure4}
\end{figure}

The expression of $\lambda_3\lambda_4+\lambda_3\lambda_5+\lambda_4\lambda_5$ and $\lambda_3\lambda_4\lambda_5$ can be determined with the parameter $m$. Let
\begin{align}
f_1(m) &=\lambda_3+\lambda_4+\lambda_5 \nonumber \\
&=\frac{1}{2}(2\sqrt{3}m^2+(9\sqrt{3}+2)m+5), \nonumber \\
f_2(m) & =\lambda_3\lambda_4+\lambda_3\lambda_5+\lambda_4\lambda_5 \nonumber\\
& =\frac{1}{2}(18m^3-18m^2+11\sqrt{3}m^2+15\sqrt{3}m+5m+3),\nonumber\\
f_3(m)&=\lambda_3\lambda_4\lambda_5 \nonumber \\
&=-\frac{27\sqrt{3}m^3}{2}+\frac{21\sqrt{3}m^2}{4}-\frac{27m^2}{4}+\frac{3m}{2}+\frac{9\sqrt{3}m}{4}+\frac{45m^3}{4}. \nonumber
\end{align}

This imples that $\lambda_1$, $\lambda_2$ and $\lambda_3$ are three
roots of the
following cubic equation:
$$\lambda^3-f_1(m)\lambda^2+f_2(m)\lambda-f_3(m)=0,$$
Solving about equation, we obtain explicitly $\lambda_3,\lambda_4,\lambda_5$
as functions of $m$. We sketch the graph of the solutions in Figure
\ref{fig:figure4}.

\

We remark that, without the symmetry considerations, finding
eigenvalues analytically, or even numerically, of an 
$8\times 8$ matrix with parameter would be very difficult, if at all
possible.

As we have pointed out before, two zero eigenvalues are expected, these
correspond to rotation and scaling invariance of the function
$\sqrt{I} U$. A central configuration is said to be {\em degenerate}\/
if there are more than two zero eigenvalues. This indeed happens when
{$f_3=\lambda_3\lambda_4\lambda_5=0$}. We conclude that for
$m^*=\frac{2\sqrt{3}+9}{18\sqrt{3}-15},$ the Hessian have one
additional zero eigenvalue. This implies that there exists a
degenerate central configuration in the $4$-body problem, as discovered
in Palmore's work \cite{Palmore1975Classifying}. See also the last
section of Xia \cite{Xia1991}.

\section{Stability of relative equilibria}
In this section, we study the local orbital structure of these
relative equilibria corresponding to the two classes of central
configurations we have discussed.

We first remark that linearization at the relative equilibria is no
longer an $8\times 8$ matrix, but instead a $16\times 16$ matrix. The
dihedral symmetry alone may not provide enough information towards
solving the linear differential equation. Additinal considerations has
to be made.

\subsection{Triangle + center problem}
The first case, the configuration with equilateral triangle plus the
center point. As
discussed in section 3.3, the $4$-body central configuration is
$$z_0=(1,0,-\frac{1}{2},\frac{\sqrt{3}}{2},- \frac{1}{2},-
\frac{\sqrt{3}}{2},0,0).$$
The corresponding relative
equilibrium solution is
$$ {\bf q}_1
=\begin{pmatrix}
cos \omega t \\
sin \omega t
\end{pmatrix},
\quad  {\bf q}_2
=\begin{pmatrix}
cos(\omega t+\frac{2 \pi}{3}) \\
sin(\omega t+\frac{2 \pi}{3})
\end{pmatrix},$$
$$  {\bf q}_3
=\begin{pmatrix}
cos(\omega t+\frac{4 \pi}{3}) \\
sin(\omega t+\frac{4 \pi}{3})
\end{pmatrix},
\quad  {\bf q}_4
=\begin{pmatrix}
0 \\
0
\end{pmatrix},$$ where $\omega^2 = \frac{\sqrt{3}}{3}+m$ according to
equation (\ref{omega}). Under uniform rotating
coordinates~\cite{1981Periodic}, the circular solution becomes
restpoint at
$$ {\bf x}_1 
=\begin{pmatrix}
1 \\
0
\end{pmatrix},
\quad {\bf x}_2
=\begin{pmatrix}
-\frac{1}{2} \\
\frac{\sqrt{3}}{2}
\end{pmatrix},$$
$$ {\bf x}_3
=\begin{pmatrix}
-\frac{1}{2} \\
-\frac{\sqrt{3}}{2}
\end{pmatrix},
\quad {\bf x}_4
=\begin{pmatrix}
0 \\
0
\end{pmatrix}.$$

To study the stability of this relative equilibrium, we need to
linearize the equation of motion around the resypoint. Recall the equation
of motion (\ref{con:inventoryflow}),

\begin{equation}
\ddot{x_i}=2\omega J\dot{x_i}+\omega^2
x_i+\frac{1}{m_i}\sum_{\substack{j=1 \\ j \ne i}
}^{n}\frac{m_im_j(x_j-x_i)}{|x_j-x_i|^3}, \nonumber
\end{equation}
for $i=1, \ldots, 4$. All terms are linear except $$\sum_{\substack{j=1 \\ j \ne i}
}^{n}\frac{m_im_j(x_j-x_i)}{|x_j-x_i|^3},$$
which is the gradient of potential function $U$. To linearize, we compute the Hessian
of $U$ at
$$z_0=(1,0,-\frac{1}{2},\frac{\sqrt{3}}{2},-\frac{1}{2},-\frac{\sqrt{3}}{2},0,0).$$
Denoting the resulting matrix as $H_3$, we have
$$\begin{pmatrix}
\frac{5\sqrt{3}}{18}+2m & 0 & -\frac{5\sqrt{3}}{36} & \frac{1}{4} & -\frac{5\sqrt{3}}{36} & -\frac{1}{4} & -2m & 0 \\
0 & -\frac{\sqrt{3}}{18}-m & \frac{1}{4} & \frac{\sqrt{3}}{36} & -\frac{1}{4} & \frac{\sqrt{3}}{36} & 0 & m \\
-\frac{5\sqrt{3}}{36} & \frac{1}{4} & \frac{\sqrt{3}}{36}-\frac{m}{4} & -\frac{1}{4}-\frac{3\sqrt{3}m}{4} & \frac{\sqrt{3}}{9} & 0 & \frac{m}{4} & \frac{3\sqrt{3}m}{4} \\
\frac{1}{4} & \frac{\sqrt{3}}{36} & -\frac{1}{4}-\frac{3\sqrt{3}m}{4} & \frac{7\sqrt{3}}{36}+\frac{5m}{4} & 0 & -\frac{2\sqrt{3}}{9} & \frac{3\sqrt{3}m}{4} & -\frac{5m}{4} \\
-\frac{5\sqrt{3}}{36} & -\frac{1}{4} & \frac{\sqrt{3}}{9} & 0 & \frac{\sqrt{3}}{36}-\frac{m}{4} & \frac{1}{4}+\frac{3\sqrt{3}m}{4} & \frac{m}{4} & -\frac{3\sqrt{3}m}{4} \\
-\frac{1}{4} & \frac{\sqrt{3}}{36} & 0 & -\frac{2\sqrt{3}}{9} & \frac{1}{4}+\frac{3\sqrt{3}m}{4} & \frac{7\sqrt{3}}{36}+\frac{5m}{4} & -\frac{3\sqrt{3}m}{4} & -\frac{5m}{4} \\
-2m & 0 & \frac{m}{4} & \frac{3\sqrt{3}m}{4} & \frac{m}{4} & -\frac{3\sqrt{3}m}{4} & \frac{3m}{2} & 0 \\
0 & m & \frac{3\sqrt{3}m}{4} & -\frac{5m}{4} & -\frac{3\sqrt{3}m}{4} & -\frac{5m}{4} & 0 & \frac{3m}{2}
\end{pmatrix} $$
We now proceed to find the eigenvalues and eigenvectors of the above
matrix. We use the same technique as before. The matrix is invariant
under an action of $S_3$. For any element $A$ in $S_3$, we still have that $H_3$ is invariant
under linear transformation $\mathscr{D}(A)$, i.e.,
$H_3\mathscr{D}(A)=\mathscr{D}(A)H_3$. Properly choosing eigenvectors
for $H_3$ as new coordinates, we have equations similar to equation
(\ref{character}).
\\
For $A=I,$
$$\lambda_1+\lambda_2+2(\lambda_3+\lambda_4+\lambda_5)=Tr(H_3)=\frac{2\sqrt{3}}{3}+6m.$$
For $A=R,$
$$\lambda_1+\lambda_2-(\lambda_3+\lambda_4+\lambda_5)=Tr(H_3\mathscr{D}(R))=\frac{\sqrt{3}}{6}-\frac{3m}{2}.$$
For $A=T$
$$\lambda_1-\lambda_2+0(\lambda_3+\lambda_4+\lambda_5)=Tr(H_3\mathscr{D}(T))=\sqrt{3}+3m.$$
Solving these equations, it gives
$$\lambda_1=\frac{2\sqrt{3}}{3}+2m, \quad \lambda_2=-\frac{\sqrt{3}}{3}-m,$$
$$\lambda_3+\lambda_4+\lambda_5=\frac{\sqrt{3}}{6}+\frac{5m}{2}.$$
By translation invariance of $U$, (or by calculating the characteristic
polynomial of $H_3,$),  there must be at least 
two zero {eigenvalues}. We
may assume
$\lambda_3=0$. 
With the determinant $(H_3+I)$, we {have} another equation
$$(\lambda_4+1)(\lambda_5+1)=\frac{3}{18}(2\sqrt{3}m-48m^2+\sqrt{3}+15m+6).$$
Then we obtain
$$\lambda_4=\frac{\sqrt{3}}{12}+\frac{5m}{4}+\frac{\sqrt{3-18\sqrt{3}m+1377m^2}}{12},$$
$$\lambda_5=\frac{\sqrt{3}}{12}+\frac{5m}{4}-\frac{\sqrt{3-18\sqrt{3}m+1377m^2}}{12}.$$

These are eigenvalues for the Hessian of $U$ at $z_0$. Now we can
easily diagonalize $H_3$ by choosing a new coordinate system.

However, we have another linear term involving $\dot{x}_i$ with
coefficient $J$, a rotation of $\pi/2$ of the tangent vectors under
the original coordinates. This term can't be diagonalized under the
new coordinate system. The question now is how the $8 \times 8$ matrix
$J_4 = Diagonal(J, J, J, J)$ acts on the eigenvectors of $H_3$. Let
$v_1, \ldots, v_4$ be eigenvectors of $H_3$ corresponding respectively
to eigenvalues $\lambda_1, \ldots, \lambda_4$. It is easy to see that
we may choose $v_2 = -Jv_1$ and $v_4 = -Jv_3$. i.e., if we group $v_1$
and $v_2$ together, and also $v_3$ and $v_4$ together, forming two
two-dimensional subspaces. The matrix $J_4$ takes the same form on
these two subspaces when we use respective eigenvectors as new
coordinates. Denoting the new coordinate systems as $y$, we have
\begin{equation}
\ddot{y_i}=2\omega J \dot{y_i}+(\omega^2I_2+\frac{1}{m_i}L_i)y_i, \label{character2}
\end{equation}
where
$$L_1=\begin{pmatrix}
\lambda_1 & 0 \\
0 & \lambda_2
\end{pmatrix}, \quad
L_2=\begin{pmatrix}
\lambda_3 & 0 \\
0 & \lambda_3
\end{pmatrix}, $$
Writing above second order equations as a first order linear systems, we have 
$$\dot{z_i}=B_iz_i, \quad i=1,2$$
where
$$B_i=\begin{pmatrix}
	0 & I_2 \\
	\omega^2 I_2+\frac{1}{m_i}L_i & 2\omega J  \\
      \end{pmatrix}.$$
We can directly solve above linear equations by finding the
eigenvalues of the matrix $B_i$. We obtain the following.

For $B_1$,
$$\lambda_1'=\lambda_2'=0, \quad \lambda_3'=\frac{\sqrt{3\sqrt{3}+9m}}{3}i, \quad \lambda_4'=-\frac{\sqrt{3\sqrt{3}+9m}}{3}i.$$
For $B_2$,
$$\lambda_5'=\lambda_6'=-\frac{\sqrt{3\sqrt{3}+9m}}{3}i, \quad \lambda_7'=\lambda_8'=\frac{\sqrt{3\sqrt{3}+9m}}{3}i.$$

For eigenspaces of $H_3$ corresponding to eigenvalues $\lambda_4$ and
$\lambda_5$, both double eigenvalues, the action of $J_4$ mixes the two
eigenspaces. Physically, this is due to the fact that the action of
$S_3$ misses the center particle. This does not happen for the square
4-body problem we will consider next. 

As for the subspaces span by the eigenspaces of $\lambda_4$ and
$\lambda_5$, it is invariant under $J_4$, one can easily write the
linear equation on this subspace. The result is a system of 8 first
order linear equations. It does not seem to admit a simple
solution. The symmetry from $S_3$ action provided an effective way to
reduce the original problem to a much simpler system.

\subsection{Square configuration}
The second case is the square configuration. As investigated in
section 3.2, the 4-body central configuration is 
$$z_0=(2,0,0,2,-2,0,0,-2).$$
The corresponding relative equilibrium solution is
$$ {\bf q}_1'
=\begin{pmatrix}
	2cos \omega' t \\
	2sin \omega' t
\end{pmatrix},
\quad  {\bf q}_2'
=\begin{pmatrix}
	2cos(\omega' t+\frac{\pi}{2}) \\
	2sin(\omega' t+\frac{\pi}{2})
\end{pmatrix},$$
$$  {\bf q}_3'
=\begin{pmatrix}
	2cos(\omega' t+\pi) \\
	2sin(\omega' t+\pi)
\end{pmatrix},
\quad  {\bf q}_4'
=\begin{pmatrix}
	2cos(\omega' t+\frac{3\pi}{2}) \\
	2sin(\omega' t+\frac{3\pi}{2})
\end{pmatrix},$$ 
with $\omega'^2 = \frac{\sqrt{2}}{16}+\frac{1}{32}$ according to
equation (\ref{omega}). Under uniform rotating
coordinates\cite{1981Periodic}, they become
$$ {\bf x}_1' 
=\begin{pmatrix}
	2 \\
	0
\end{pmatrix},
\quad {\bf x}_2'
=\begin{pmatrix}
	0 \\
	2
\end{pmatrix},$$
$$ {\bf x}_3'
=\begin{pmatrix}
	-2 \\
	0
\end{pmatrix},
\quad {\bf x}_4'
=\begin{pmatrix}
	0 \\
	-2
\end{pmatrix}.$$
Similarly, we linearize the equation of motion
(\ref{con:inventoryflow}) around the orbit. We compute the Hessian of
$U$ at  
$$z_0=(2,0,0,2,-2,0,0,-2).$$
Denoting the resulting Hessian as $H_4$, we have
$$
\begin{pmatrix}
	\frac{\sqrt{2}}{32}+\frac{1}{32} & 0 & -\frac{\sqrt{2}}{64} & \frac{3\sqrt{2}}{64} & -\frac{1}{32} & 0 & -\frac{\sqrt{2}}{64} & -\frac{3\sqrt{2}}{64} \\
	0 & \frac{\sqrt{2}}{32}-\frac{1}{64} & \frac{3\sqrt{2}}{64} & -\frac{\sqrt{2}}{64} & 0 & \frac{1}{64} & -\frac{3\sqrt{2}}{64} & -\frac{\sqrt{2}}{64} \\
	-\frac{\sqrt{2}}{64} & \frac{3\sqrt{2}}{64} & \frac{\sqrt{2}}{32}-\frac{1}{64} & 0 & -\frac{\sqrt{2}}{64} & -\frac{3\sqrt{2}}{64} & \frac{1}{64} & 0 \\
	\frac{3\sqrt{2}}{64} & -\frac{\sqrt{2}}{64} & 0 & \frac{\sqrt{2}}{32}+\frac{1}{32} & -\frac{3\sqrt{2}}{64} & -\frac{\sqrt{2}}{64} & 0 & -\frac{1}{32} \\
	-\frac{1}{32} & 0 & -\frac{\sqrt{2}}{64} & -\frac{3\sqrt{2}}{64} & \frac{\sqrt{2}}{32}+\frac{1}{32} & 0 & -\frac{\sqrt{2}}{64} & \frac{3\sqrt{2}}{64} \\
	0 & \frac{1}{64} &  -\frac{3\sqrt{2}}{64} & -\frac{\sqrt{2}}{64} & 0 & \frac{\sqrt{2}}{32}-\frac{1}{64} & \frac{3\sqrt{2}}{64} & -\frac{\sqrt{2}}{64} \\
	-\frac{\sqrt{2}}{64} & -\frac{3\sqrt{2}}{64} & \frac{1}{64} & 0 & -\frac{\sqrt{2}}{64} & \frac{3\sqrt{2}}{64} & \frac{\sqrt{2}}{32}-\frac{1}{64} & 0 \\
	-\frac{3\sqrt{2}}{64} & -\frac{\sqrt{2}}{64} & 0 & -\frac{1}{32} & \frac{3\sqrt{2}}{64} & -\frac{\sqrt{2}}{64} & 0 & \frac{\sqrt{2}}{32}+\frac{1}{32}
\end{pmatrix}
$$
For any elemnt $A$ in the group $D_4,$ we have $H_4$ is invariant
under linear transformation $\mathscr{D}(A),$ i.e., $H_4
\mathscr{D}(A)=\mathscr{D}(A)H_4.$. Similar to equation (\ref{trace
  d4}), we have the following equations for eigenvalues of $H_4$.
\\
For $A=e,$
$$\lambda_1+\lambda_2+\lambda_3+\lambda_4+ 2(\lambda_5+ \lambda_6) =
\frac{\sqrt{2}}{4}+\frac{1}{16}.$$
For $A=a,$
$$ \lambda_1+\lambda_2- \lambda_3-\lambda_4 =0.$$
For $A=a^2,$
$$ \lambda_1 +\lambda_2+\lambda_3+\lambda_4-2(\lambda_5+\lambda_6) = \frac{1}{16}. $$
For $A=r,$
$$ \lambda_1-\lambda_2+\lambda_3-\lambda_4 =\frac{3}{16}.$$
For $A=ar,$
$$ \lambda_1-\lambda_2-\lambda_3+\lambda_4 =-\frac{3\sqrt{2}}{8}.$$
Solving these equations, we have
$$ \lambda_1 = \frac{1}{16} + \frac{\sqrt{2}}{8}, \quad \lambda_2= -\frac{1}{32}-\frac{\sqrt{2}}{16},$$
$$ \lambda_3 = \frac{1}{16} - \frac{\sqrt{2}}{16}, \quad \lambda_4 = -\frac{1}{32}+\frac{\sqrt{2}}{8} \quad \lambda_5 +\lambda_6 = \frac{\sqrt{2}}{16}.$$
By translation invariance, or the character polynomial of $H_4,$ there
are two zero {eigenvalues}. Assuming $\lambda_5 = 0,$
then $\lambda_6= \frac{\sqrt{2}}{16}.$

Again, we need to analyze the action of $J$ on the eigenvectors of
$H_4$. This case turn out to be much simpler. The matrix $J_4$ keeps
invariant of the following four two-dimensional subspaces: the span of
the first and the second eigenvectors; the span of the third and the
fourth eigenvectors; the eigenspaces of $\lambda_5$; the eigenspace
$\lambda_6$. Using $y$ as new coordinate system for each subspaces,
the linear equation takes the following form,
(\ref{con:inventoryflow}) at $z_0$ and choose eigenvectors of $H_4$ as
a basis in $R^8$, we have
\begin{equation*}
	\ddot{y_i}=2\omega J \dot{y_i}+(\omega^2I_2+L_i')y_i,
\end{equation*}
where
$$L_1'=\begin{pmatrix}
	\lambda_1 & 0 \\
	0 & \lambda_2
\end{pmatrix}, \quad
L_2'=\begin{pmatrix}
	\lambda_3 & 0 \\
	0 & \lambda_4
\end{pmatrix}, $$
$$L_3'=\begin{pmatrix}
	\lambda_5 & 0 \\
	0 & \lambda_5
\end{pmatrix}, \quad
L_4'=\begin{pmatrix}
	\lambda_6 & 0 \\
	0 & \lambda_6
\end{pmatrix}.$$
It is easy to calculate eigenvalues for the second-order differential equations.
\\
For $i= 1,$
$$\lambda_1'=\lambda_2'=0, \quad \lambda_3'=\frac{\sqrt{2+4\sqrt{2}}}{8}i, \quad \lambda_4'= -\frac{\sqrt{2+4\sqrt{2}}}{8}i. $$
For $i=2,$
$$\lambda_5'= \frac{\sqrt{\sqrt{68\sqrt{2}-9}i-2\sqrt{2}-1}}{8}, \quad \lambda_6'= -\frac{\sqrt{\sqrt{68\sqrt{2}-9}i-2\sqrt{2}-1}}{8},$$
$$\lambda_7'= \frac{\sqrt{-\sqrt{68\sqrt{2}-9}i-2\sqrt{2}-1}}{8}, \quad \lambda_8'= -\frac{\sqrt{-\sqrt{68\sqrt{2}-9}i-2\sqrt{2}-1}}{8}.$$
For $i=3,$
$$\lambda_9'= \lambda_{10}'=\frac{\sqrt{2+4\sqrt{2}}}{8}i, \quad \lambda_{11}'= \lambda_{12}'=-\frac{\sqrt{2+4\sqrt{2}}}{8}i.$$
For $i=4,$
$$\lambda_{13}'= -\frac{2^{\frac{1}{4}}}{4}+\frac{\sqrt{2+4\sqrt{2}}}{8}i, \quad \lambda_{14}'=\frac{2^{\frac{1}{4}}}{4}+\frac{\sqrt{2+4\sqrt{2}}}{8}i, $$
$$\lambda_{15}'= -\frac{2^{\frac{1}{4}}}{4}-\frac{\sqrt{2+4\sqrt{2}}}{8}i, \quad \lambda_{16}'=\frac{2^{\frac{1}{4}}}{4}-\frac{\sqrt{2+4\sqrt{2}}}{8}i.$$

Therefore,
$\lambda_5',\lambda_6',\lambda_7',\lambda_8,\lambda_{13}',\lambda_{14}',\lambda_{15}',\lambda_{16}'$
all have nonzero real parts, we conclude that the square configuration is unstable. Eigenvalues $
\lambda_3',\lambda_4', \lambda_9', \lambda_{10}', \lambda_{11}',
\lambda_{12}'$ are pure imaginary numbers. This implies that the
elliptic direction has dimension $6$.

\bibliographystyle{acm}

\begin{thebibliography}{10}

\bibitem{2012Finiteness}
{\sc Albouy, A., and Kaloshin, V.}
\newblock Finiteness of central configurations of five bodies in the plane.
\newblock {\em Ann. of Math. (2) 176}, 1 (2012), 535--588.

\bibitem{2012Four}
{\sc Cors, J.~M., and Roberts, G.~E.}
\newblock Four-body co-circular central configurations.
\newblock {\em Nonlinearity 25}, 2 (2012), 343--370.

\bibitem{2019On}
{\sc Fernandes, A.~C., Mello, L.~F., and Vidal, C.}
\newblock On the uniqueness of the isoceles trapezoidal central configuration
  in the 4-body problem for power-law potentials.
\newblock {\em Nonlinearity 33}, 1 (2020), 388--407.

\bibitem{HamptonMoeckel2006}
{\sc Hampton, M., and Moeckel, R.}
\newblock Finiteness of relative equilibria of the four-body problem.
\newblock {\em Invent. Math. 163}, 2 (2006), 289--312.

\bibitem{leandro2017factorization}
{\sc Leandro, E.~S.}
\newblock Factorization of the stability polynomials of ring systems.
\newblock {\em arXiv preprint arXiv:1705.02701\/} (2017).

\bibitem{MR3951830}
{\sc Leandro, E. S.~G.}
\newblock Structure and stability of the rhombus family of relative equilibria
  under general homogeneous forces.
\newblock {\em J. Dynam. Differential Equations 31}, 2 (2019), 933--958.

\bibitem{1981Periodic}
{\sc Meyer, K.~R.}
\newblock Periodic solutions of the {$N$}-body problem.
\newblock {\em J. Differential Equations 39}, 1 (1981), 2--38.

\bibitem{easton1993introduction}
{\sc Meyer, K.~R., Hall, G.~R., and Offin, D.}
\newblock {\em Introduction to {H}amiltonian dynamical systems and the
  {$N$}-body problem}, second~ed., vol.~90 of {\em Applied Mathematical
  Sciences}.
\newblock Springer, New York, 2009.

\bibitem{MR1350320}
{\sc Moeckel, R.}
\newblock Linear stability analysis of some symmetrical classes of relative
  equilibria.
\newblock In {\em Hamiltonian dynamical systems ({C}incinnati, {OH}, 1992)},
  vol.~63 of {\em IMA Vol. Math. Appl.} Springer, New York, 1995, pp.~291--317.

\bibitem{2017On}
{\sc Oliveira, A., and Cabral, H.}
\newblock On stacked central configurations of the planar coorbital satellites
  problem.
\newblock {\em Discrete Contin. Dyn. Syst. 32}, 10 (2012), 3715--3732.

\bibitem{Palmore1975Classifying}
{\sc Palmore, J.~I.}
\newblock Classifying relative equilibria. {II}.
\newblock {\em Bull. Amer. Math. Soc. 81\/} (1975), 489--491.

\bibitem{Palmore1976}
{\sc Palmore, J.~I.}
\newblock Measure of degenerate relative equilibria. {I}.
\newblock {\em Ann. of Math. (2) 104}, 3 (1976), 421--429.

\bibitem{MR1816912}
{\sc Roberts, G.~E.}
\newblock Linear stability in the {$1+n$}-gon relative equilibrium.
\newblock In {\em Hamiltonian systems and celestial mechanics ({P}\'{a}tzcuaro,
  1998)}, vol.~6 of {\em World Sci. Monogr. Ser. Math.} World Sci. Publ., River
  Edge, NJ, 2000, pp.~303--330.

\bibitem{Saari1980On}
{\sc Saari, D.~G.}
\newblock On the role and the properties of {$n$}-body central configurations.
\newblock {\em Celestial Mech. 21}, 1 (1980), 9--20.

\bibitem{smale1998mathematical}
{\sc Smale, S.}
\newblock Mathematical problems for the next century.
\newblock {\em Math. Intelligencer 20}, 2 (1998), 7--15.

\bibitem{2012Representation}
{\sc Steinberg, B.}
\newblock {\em Representation theory of finite groups}.
\newblock Universitext. Springer, New York, 2012.
\newblock An introductory approach.

\bibitem{Xia1991}
{\sc Xia, Z.}
\newblock Central configurations with many small masses.
\newblock {\em J. Differential Equations 91}, 1 (1991), 168--179.

\bibitem{xia2004convex}
{\sc Xia, Z.}
\newblock Convex central configurations for the {$n$}-body problem.
\newblock {\em J. Differential Equations 200}, 2 (2004), 185--190.

\bibitem{Xia2008}
{\sc Xia, Z.}
\newblock {Symmetries in N-body problem}.
\newblock {\em AIP Conference Proceedings 1043}, 126 (2008), 126--132.

\end{thebibliography}

\end{document}